\input amstex
\input amsppt.sty
\magnification=\magstep1
\overfullrule = 0pt
\TagsOnRight
\NoBlackBoxes
\leftheadtext{V\. Pestov}
\rightheadtext{Universal constructions}
\def\norm #1{{\left\Vert\,#1\,\right\Vert}}
\def\e {{\varepsilon}}
\def\R {{\Bbb R}}
\def\C {{\Bbb C}}
\def\K {{\Bbb K}}
\def\Q {{\Bbb Q}}
\def\Z{{\Bbb Z}}
\def\N{{\Bbb N}}
\def\Homeo{{\operatorname{Homeo}\,}}
\def\St{{\operatorname{St}}}
\def\Aut{{\operatorname{Aut}\,}}
\def\Sym{{\operatorname{Sym}\,}}
\def\Iso{{\operatorname{Iso}\,}}
\def\Spec{{\operatorname{Spec}\,}}
\def\dist{{\operatorname{dist}\,}}
\def\supp{{\operatorname{supp}\,}}
\def\tri
{\ifhmode\unskip\nobreak\fi\ifmmode\ifinner\else\hskip 5pt\fi\fi
 \hfill\hbox{\hskip 5pt$\blacktriangle$\hskip 1pt}}
\def\quest
{\ifhmode\unskip\nobreak\fi\ifmmode\ifinner\else\hskip 5pt\fi\fi
 \hfill\hbox{\hskip 5pt$\boxed{\text{\sixbf ?}}$\hskip 1pt}}
\def\newqed
{\ifhmode\unskip\nobreak\fi\ifmmode\ifinner\else\hskip 5pt\fi\fi
 \hfill\hbox{\hskip 5pt$\square$\hskip 1pt}}

\rightline{\eightit Research Report RP-97-219 (April 1997),}
\rightline{\eightit School of Mathematical and Computing Sciences,}
\rightline{\eightit Victoria University of Wellington}
\vskip 1cm
\topmatter 
\title
Some universal constructions in abstract topological dynamics
\endtitle 
\author
Vladimir Pestov
\endauthor
\address School of Mathematical and Computing Sciences, 
Victoria University of Wellington,
P.O. Box 600, Wellington, New Zealand 
\endaddress 
\email vladimir.pestov$\@$mcs.vuw.ac.nz
\endemail
\affil Victoria University of Wellington, 
Wellington, New Zealand
\endaffil
\dedicatory Dedicated to Professor Robert Ellis on the occasion of his
retirement
\enddedicatory
 
\abstract This small survey of basic universal constructions
related to the actions of topological groups on compacta is centred around
a new result  --- an intrinsic description of extremely amenable
topological groups (i.e., those having a fixed point in each compactum
they act upon), 
solving a 1967 problem by Granirer. Another old problem whose solution
(in the negative) is noted here, is a 1957 conjecture by Teleman
on irreducible representations of general topological groups.
Our exposition covers the greatest ambit
and the universal minimal flow, as well as closely related
constructions and results from representation
theory. 
We present in a simplified fashion the
known examples of extremely amenable groups and discuss
their relationship with a 1969 problem by Ellis.
Open questions, including those from 
bordering disciplines, are reviewed along the way.
\endabstract
\endsubjclass

\endtopmatter 
\document 
\head 0. Introduction
\endhead
The principal object of study in abstract topological dynamics is a
$G$-flow (topological transformation group, topological dynamical
system), that is, a triple $\frak X= (G,X,\tau)$
consisting of a Hausdorff topological group $G$, a Tychonoff
topological space $X$, and a continuous action $\tau\colon G\times X\to X$.
Theory is at its most advanced for compact phase spaces $X$,
in which case $\tau$ can be thought of as a continuous
homomorphism $G\to\Homeo X$, where the group of self-homeomorphisms is equipped
with the compact-open topology. It is therefore most natural to 
begin a journey into the realm of abstract topological dynamics
by asking: does {\it every} topological group
$G$ admit an effective action on a compact space, thereby allowing for a
dynamically substantial theory?

This fundamental question was answered in the positive forty years
ago by Teleman 
\cite{T} who produced a fine collection of results
on the crossroads of topological dynamics with representation
theory and theory of $C^\ast$-algebras. We present those results,
leading to the concept of the greatest ambit, in Sections  1 and 4-5.
Among the later corollaries, is a
1986 Uspenski\u\i's observation \cite{U} of the existence of a
universal Polish group, answering the `right' version of Schreier--Ulam's
question from the Scottish Book (Section 3). We also point out, possibly for
the first time, that a 1957 conjecture of Teleman can be refuted through 
combining some
developments in operator theory and in abstract harmonic analysis
(Section 2).

Of all compact $G$-flows, minimal $G$-flows (that is, as name suggests,
those having no proper compact $G$-subflows) are of particular
importance. In Section 6 we introduce the universal minimal flow, 
$\Cal M(G)$, for any topological group $G$. The next important 
question to address
is that of the `size' of the universal minimal flow. In
particular, when is
$\Cal M(G)$ nontrivial, and when is the action of $G$ upon $\Cal M(G)$ 
effective? 

A topological group $G$ such that the universal (therefore, every)
minimal $G$-flow is a singleton
is called {\it extremely amenable} \cite{Gr1, Mi, Gr2}
(or else a group having the {\it fixed point on compacta property} 
\cite{G3}).

The very existence of such groups was at first an open problem, 
though topological semigroups
with a similar property were known to exist in abundance 
\cite{Gr1, Mi}. By now there are known examples of extremely amenable
topological groups of at least three different kinds
\cite{HR, G3, P}. 
In 1967 Granirer \cite{Gr1} suggested 
the problem of describing
extremely amenable topological groups in intrinsic terms. 
We propose such a description in Sections 7 and 8. 
Recall that a subset $S$ of a topological group $G$ is
{\it left syndetic} if $KS=G$ for some compact $K\subseteq G$.
We show that a topological group $G$ is extremely amenable if and only
if for every left syndetic subset $S\subseteq G$ the set
$SS^{-1}$ is everywhere dense in $G$. (One can also replace 
`left syndetic' with `big on the left,' that is, left syndetic in the
discrete topology.)

Employing this new criterion, we sketch the
known examples of extremely amen\-able groups. Historically the first
such example, constructed by Herer and Christensen,
utilises the existence of abelian topological groups without nontrivial
unitary representations (see Section 13). 
Another source of examples, presented in Section 9,
was discovered by Glasner and (independently) Furstenberg
and Weiss, who built
up on the work of Gromov and Milman, and it includes the (even monothetic)
group of measurable maps from $I$ to the unit circle, equipped with the
$L_1$-metric. Finally, among the extremely amenable groups proposed by the
present author (Section 10), is the group $\Homeo_+(I)$ of 
orientation-preserving homeomorphisms
of the unit interval. We review a number of dynamical 
consequences for groups
having importance in Analysis and observe that beyond the locally compact
case, amenability is not inherited by closed subgroups (Section
10).

The existence of extremely amenable groups answers in the negative a
1969 question by Ellis from \cite{E3}, which was 
largely responsible
for research in the area and can be reformulated in equivalent terms as
follows: are points of the greatest ambit $\Cal S(G)$ 
of a topological group $G$ separated by continuous equivariant maps
to minimal $G$-flows? This problem and related developments form
the subject of Section 11.

Our criterion of extreme amenability is accompanied by a satellite
result describing, also for the first time,
those topological groups $G$ admitting effective minimal actions on
compacta (Section 12). 
A topological group $G$ acts effectively on $\Cal M(G)$ if and only
if for every $g\in G$, $g\neq e$ there is a left syndetic set
$S$ with $x\notin \overline{SS^{-1}}$. 
We derive from this criterion an
easy proof of the fact that every locally compact group
acts effectively on the universal minimal flow. This is of course
a corollary of
a stronger theorem by Veech \cite{V2}
stating that every locally
compact group acts {\it freely} on a compactum --- 
the result of which a `soft' proof has so far escaped us.

Finally, in Section 13 we discuss a few results and questions related to the
old and still open problem of describing those topological
groups admitting a complete system of strongly continuous
unitary representations in Hilbert spaces.

Our usage of basic concepts from topology, functional analysis,
topological algebra, representation theory,
and topological dynamics is fairly standard,
while our list of bibliography is kept to a bare minimum.

\subhead Acknowledgements  \endsubhead
The present article has grown out of author's presentations of 
some of the results at
the University of Maryland (September 1996), at the Semi-Annual Workshop on
Dynamical Systems, Penn State University (October 1996), and
at the Tsukuba Topological Conference, Tsukuba Science City, Japan
(December 1996).
The author is most thankul to Professors Joe Auslander and
Dmitri Shakhmatov, as well as to the Organizers of the above Conferences,
for making those presentations possible. Most stimulating e-mail discussions 
with Prof. Eli Glasner and important bibliographical 
remarks by Prof. Alexander Kechris and Dr. Michael Megrelishvili
are also acknowledged. Finally, a number of useful
improvements were suggested by the referee of this article.

\head  1. Flows and Representations \endhead

Topological dynamics is intimately linked to representation theory,
as can be seen from the following observation.
Let $\rho$ be a representation of an abstract group $G$ in a Banach
space $E$ by isometries, that is, $\rho\colon G\to\Iso (E)$ is an
homomorphism to the group of all linear isometries of $E$ onto itself.
For every $g\in G$ the dual operator, $\rho_g^\ast\colon E^\ast\to E^\ast$,
to the isometry $\rho_g\colon E\to E$ is an isometry again, and
the restriction, $\rho_g^B$, of $\rho^\ast_{g^{-1}}$ to the unit ball,
$B^\ast_w$, of $E^\ast$, equipped with the $w^\ast$-topology
(the coarsest topology making each evaluation map
$\hat x\colon E^\ast\to\K$, $\hat x(f)=f(x)$, continuous), is a
self-homeomorphism. It is easy to check that the resulting map
$$\rho^B\colon G\ \ni g\mapsto \rho_g^B\in\Homeo(B^\ast_w)$$
is a group homomorphism, that is, an action of $G$ on the compact space
$B^\ast_w$. 

Now assume that $G$ is a topological group. Recall that a representation
$\rho\colon G\to\operatorname{GL}(E)$
(in particular, a representation $\rho\colon G\to\Iso(E)$)
is {\it strongly continuous}
if it is continuous with respect to the topology on
$\operatorname{GL}(E)$ inherited from the
Tychonoff product topology on $E^E$. (This, the so-called
{\it strong operator topology,} though not a group topology on 
$\operatorname{GL}(E)$, becomes a group topology when restricted to 
$\Iso(E)$.)
\proclaim{PROPOSITION  1.1 {\rm (Teleman \cite{T})}} {\rm (i)} If
$\rho$ is strongly continuous then $\rho^B$ is continuous.
\item{\rm (ii)} If $\rho$ is faithful then $\rho^B$ is effective.
\item{\rm (iii)} If $\rho\colon G\hookrightarrow\Iso (E)$ is an embedding of
topological groups, then so is 
$\rho^B\colon G\hookrightarrow\Homeo(B^\ast_w)$.
\endproclaim
\demo{PROOF} Let us compare 
standard basic neighbourhoods of identity in the topological group
$\Iso(E)$,
$$\{ u\in\Iso(E)\colon\forall i=1,2,\dots,n,~
\norm{u(x_i)-x_i}<\e\},~~x_i\in E,~\e>0,$$
and in $\Homeo(B^\ast_w)$,
$$\{u\in\Homeo(B^\ast_w)\colon \forall f\in B^\ast_w,~~
(u(f),f)\in V\},~~ V\in\Cal U_X,$$
where $\Cal U_X$ is the (unique) compatible uniformity on 
$B^\ast_w$. Applying the Hahn--Banach theorem in the first case
and noticing that $\Cal U_X$ is 
induced by the additive uniformity on
$E^\ast$ equipped with the weak$^\ast$ topology in the second case,
we reduce both neighbourhoods to the same form: 
$$\{u \colon \forall i=1,2,\dots, n,~\forall f\in B^\ast_w,~~
\vert f(u(x_i)) - f(x_i)\vert <\e\},~ x_i\in E, ~\e>0 $$
It implies that $\Iso(E)$ canonically embeds into $\Homeo(B^\ast_w)$ as
a topological subgroup, whence all three statements follow. \newqed\enddemo
\smallskip

To conclude that every topological group $G$ acts effectively on
a compact space, it is enough to observe the following.

\proclaim{THEOREM 1.2 {\rm (Teleman \cite{T})}} 
Every topological group $G$ admits a faithful strongly continuous
representation $\rho$ by isometries in a suitable Banach space $E$.
Moreover, $\rho\colon G\to\Iso(E)$ is in this case a topological
embedding.
\endproclaim

\demo{PROOF} We avoid using the terms `left' and
`right' uniform structure altogether because there is a
considerable confusion surrounding their usage.
Instead we denote, once and for all, by
$\Cal U_\Rsh(G)$ the uniform structure on $G$ whose base is formed by
the entourages of the diagonal
$$V_\Rsh =\{(x,y)\colon xy^{-1}\in V\},$$
where $V$ is a neighbourhood
of $e_G$. We say that a function $f\colon G\to\C$ is
$\Cal U_\Rsh${\it -uniformly continuous,} or 
$\Cal U_\Rsh${\it -u.c.,} if 
$$\forall \e>0,~\exists V\in\Cal N(G),~\forall x,y\in G,~ 
xy^{-1}\in V\Rightarrow
\vert f(x)-f(y)\vert<\e \tag{$\ast$}$$
The uniformity $\Cal U_\Lsh$ and $\Cal U_\Lsh$-uniformly continuous
functions on a topological group are introduced in a similar way,
though they are never used in this article.

The totality of all bounded complex-valued $\Cal U_\Rsh$-u.c. functions
on a topological group $G$ equipped with the supremum norm and
pointwise operations is easily checked to form a commutative
$C^\ast$-algebra --- indeed a $C^\ast$-subalgebra of
$C(\beta G)\cong l_\infty(\vert G\vert)$ --- which we denote by
$C^\ast_\Rsh(G)$.

Every element $g\in G$ determines a left shift $L_g\colon C^\ast_\Rsh(G)
\to C^\ast_\Rsh(G)$ via
$$L_g(f)(h)=f(g^{-1}h)$$
The emerging representation $L\colon G\to\Iso C^\ast_\Rsh(G)$ is
strongly continuous: if $V$ is symmetric and as in ($\ast$), 
then $\norm {L_gf-f}<\e$ whenever $g\in V$. Moreover,
$L$ is a topological embedding: given a symmetric
neighbourhood of
identity, $V$, one can construct a $\Cal U_\Rsh$-u.c. bounded
$f\colon G\to\C$ with $f\equiv 0$ outside $V$ and
$f(e)=1$; the set $W=\{u\in \Iso C^\ast_\Rsh(G)\colon 
\norm{u(f)-f}<1\}$ is open in the strong topology and
$e\in W\cap G\subseteq V$. \newqed\enddemo
\smallskip

It seems that Teleman's results later became a part of
the folklore and were even sometimes reproduced with proofs
without referring to the original. 

\head  2. Teleman's Conjecture on Irreducible Banach Representations\endhead
Having in view the celebrated
1943 Gelfand--Ra\u\i kov theorem on the existence of a complete system 
of irreducible unitary representations for any locally compact group,
Teleman asked in his 1957 paper \cite{T} whether every topological
group admits a complete system of strongly continuous irreducible
representations in Banach spaces by isometries. A negative answer to
this question can be deduced from more recent results obtained
in operator theory and abstract harmonic analysis, but it is not apparent
if it was ever noticed by anyone.

Recall that a topological group is called {\it minimally almost periodic}
if it admits no nontrivial finite-dimensional unitary representations,
and {\it monothetic} if it contains an everywhere dense cyclic subgroup.

\proclaim{PROPOSITION 2.1} A monothetic minimally almost periodic 
topological group $G$ admits no nontrivial irreducible strongly
continuous representations in a complex Banach space by isometries.
\endproclaim

\demo{PROOF} Assume the existence of such a representation, $\rho$.
Let $x$ be a generator of an everywhere dense cyclic subgroup of $G$.
The operator $\rho_x$ is then an isometry of the complex Banach space
of representation, $E$, and therefore, as a well-known consequence
of Riesz' theorem \cite{B, th. 3.1}, admits a
nontrivial closed invariant subspace, $F$, whenever $\dim E>1$.
Since $F$ remains invariant for all operators
$\rho_g$, $g\in G$, one concludes that $\dim E=1$ and the representation
$\rho$, being a continuous character $G\to U(1)$, is trivial.
\newqed\enddemo
\smallskip
Since there are numerous known examples of minimally almost periodic 
monothetic topological groups (see e.g.
\cite{DPS}), Teleman's conjecture does not survive.
However, it can be brought back to life through simply dropping the
requirement `by isometries.'

\definition{QUESTIONS  2.2} (i) Does every topological group admit a complete
system of strongly continuous irreducible representations in complex
Banach spaces? 

(ii) Assuming a topological group $G$ admits a faithful strongly
continuous representation in a Hilbert space, does it admit a complete 
system of irreducible representations in Hilbert spaces?
\quest
\enddefinition

If answered positively, the second question, applied to the
topological group from Glasner--Furstenberg--Weiss's example 
(see Section 9 below), would provide
a solution to the Invariant Subspace Problem.

\head  3. Universal Group of Countable Weight\endhead
It was asked by Schreier and
Ulam in the Scottish Book (problem 103 in the 1981 
edition \cite{Ma}) whether or
not there existed a universal separable topological group, that is, a
separable topological group $G$ such that every other such group is
isomorphic to a topological subgroup of $G$. The answer is `no' because
of simple set-theoretic reasons: while a separable group can have at
most 
$\frak c$ pairwise distinct subgroups, there are
$\exp\frak c$ pairwise non-isomorphic separable topological groups,
cf. the comments on pp. 184--186 in \cite{Ma}. 
However, the following elegant result seems to answer 
in the positive the
`right' version of Schreier--Ulam's question as suggested by 
Kallman ({\it ibid.}). 

\proclaim{THEOREM  3.1 {\rm (Uspenski\u\i\  \cite{U})}}
The group $\Homeo(I^{\aleph_0})$ is a universal separable metrizable
topological group. In other words, every separable metrizable
topological group is isomorphic with a topological subgroup of
$\Homeo(I^{\aleph_0})$.
\endproclaim

\demo{PROOF} If $G$ is separable metrizable, then its points 
are separated from closed subsets 
by elements of a countable subset --- and therefore by those of a
separable Banach subspace, say $E$ --- of
$C^\ast_\Rsh(G)$, and such a subspace can be also assumed
invariant under left translations. The action of $G$ on $E$ determines
an isometric embedding $G\hookrightarrow\Iso(E)$, and Proposition 1.1
yields an embedding of topological groups $G\hookrightarrow\Homeo(B^\ast_w)$.
How it remains to apply Keller's theorem \cite{BP}: all convex compact
infinite-dimensional metrizable subsets of locally convex spaces are
homeomorphic to the Hilbert cube $Q=I^{\aleph_0}$. 
\newqed\enddemo

\head  4. Dynamics and $C^\ast$-Algebras\endhead
Let us cast a closer look at Teleman's construction. So far we made
no use of the fact that $C^\ast_\Rsh(G)$ is a commutative 
unital $C^\ast$-algebra
rather than just a Banach space, and the shifts $L_g$ are
$C^\ast$-algebra automorphisms and not merely isometries. 

Assume for a while that a topological group $G$ acts strongly continuously
on a commutative unital $C^\ast$-algebra $A$ by automorphisms.
The Gelfand space, $\Spec_\C A$, formed by all multiplicative linear 
functionals $f\colon A\to\C$, is a closed subspace of the unit 
ball $B^\ast_w$, and it is clearly invariant under the action 
$\rho^B$ of $G$. Therefore, $G$ acts on $\Spec_\C A$
continuously. Going in the opposite direction, if $X$ is a compact
$G$-flow, then $G$ acts strongly continuously by left shifts on the
commutative $C^\ast$-algebra $C(X)$ of all complex-valued continuous
functions on $X$, and $X\cong\Spec_\C C(X)$ as $G$-spaces.
The functorial correspondence $A\mapsto\Spec_\C A$ is an
anti-equivalence of categories.

\proclaim{THEOREM 4.1} Let $G$ be a topological group. The category of
all compact $G$-flows and $G$-equivariant continuous maps
is anti-equivalent to
the category of all strongly continuous
representations of $G$ in commutative unital $C^\ast$-algebras by
automorphisms and
intertwining unital $C^\ast$-algebra morphisms.
\newqed\endproclaim

This simple observation helps to explain why
theory of strongly continuous actions of 
topological groups on
{\it arbitrary} unital $C^\ast$-algebras by automorphisms is universally
accepted as the adequate version of
{\it non-commutative dynamics}, cf. \cite{Gu}.

\head  5. The Greatest Ambit\endhead

For a topological group $G$, denote the Gelfand space 
$\Spec_\C C_\Rsh(G)$ by $\Cal S(G)$. Each element
$g\in G$ determines an evaluation functional, $\hat g\in\Cal S(G)$,
via $\hat g(f)=f(g)$.
Since $\Cal U_\Rsh$-u.c. functions separate points and closed subsets
in a topological group, the resulting canonical map 
$G\ni g\mapsto \hat g\in\Cal S(G)$ is a homeomorphic embedding. It is
also equivariant with respect to the action of $G$ on itself by
left translations: $L^\ast_{g^{-1}}(\hat x) =\widehat{gx}$ for all
$g,x\in G$. For this reason, we will be identifying each element
$g\in G$ with $\hat g\in\Cal S(G)$.
Notice that the $G$-flow $\Cal S(G)$ has a distinguished
point, $e$. The $G$-orbit of $e$ is everywhere dense in $\Cal S(G)$,
since the latter space is a compactification of $G$; it is, in fact,
the Samuel compactification (see e.g. \cite{En}) 
of the uniform space $(G,\Cal U_\Rsh(G))$. 
\definition{DEFINITION 5.1} A compact $G$-flow with a distinguished point,
$(X,x)$, is called a $G$-{\it ambit} if the orbit of $x$ is
everywhere dense in $X$. A {\it morphism} between $G$-ambits
$(X,x)$ and $(Y,y)$ is a continuous $G$-equivariant map 
$X\to Y$ taking $x$ to $y$.
\tri\enddefinition

Let $X=(X,x)$ be an arbitrary $G$-ambit for a topological group $G$.
Every $f\in C(X)$ determines a function $\tilde f\colon G\to\C$ 
via `pullback:' $\tilde f(g)=f(gx)$. It is easy to see that $\tilde f$
is bounded and $\Cal U_\Rsh$-u.c.
The $C^\ast$-algebra monomorphism
$C(X)\hookrightarrow C^\ast_\Rsh(G)$, determined by the
correspondence $f\mapsto\tilde f$, intertwines the representations of
$G$ in $C(X)$ and $C^\ast_\Rsh(G)$ and therefore
determines, according to our earlier
observation, a dual morphism of $G$-flows
$\Cal S(G)\to X$, which takes $e$ to $x$ and is therefore onto.
(The image of $\Cal S(G)$ is compact and therefore closed, and is at the
same time everywhere dense since it contains the orbit
$G\cdot x$.) We come to the following.

\proclaim{THEOREM 5.2 {\rm (Teleman)}} Every topological group $G$ possesses
the greatest ambit $\Cal S(G)=(\Cal S(G),e)$ in the sense that every
$G$-ambit is the image of $\Cal S(G)$ under a unique morphism of
$G$-ambits. The greatest ambit contains $G$ itself as an
everywhere dense $G$-subspace.
\newqed\endproclaim

As we will see later, the concept is very useful. It would be interesting
to know the following.

\definition{QUESTION  5.3} Is there an analogue of 
the greatest ambit in noncommutative dynamics?
\quest\enddefinition

\example{EXAMPLE 5.4}
Let us describe two `diametrically opposite' situations where the
structure of the greatest ambit is relatively well understood.

(i) If $G$ is precompact, then $\Cal S(G)$ is simply the compact topological
group, $\widehat G$, completion of $G$, equipped with the left action of
the subgroup $G$. 

(ii) If $G$ is discrete, then every complex-valued function
on $G$ is $\Cal U_\Rsh$-u.c. and as a corollary, $\Cal S(G)$ is
$\beta G$, the Stone-\u Cech compactification of $G$, equipped with
the natural action of $G$ on elements of $\beta G$ (which can be
identified with
all ultrafilters on the set $G$) by left translations: if $g\in G$
and $U\in\beta G$, then $gU=\{gA\colon A\in U\}\in\beta G$. 
(See \cite {A, E2, Vr}.)
\tri\endexample

The situation `in between' is much less thoroughly understood.
One significant recent advance is the description by Turek
\cite{Tu} of the greatest ambit, $\Cal S(\R)$,
of the additive group of reals, $\R$, with the usual topology.

It only takes a slight modification of the arguments used in 1.1 and 1.2 
to establish the following.

\proclaim{THEOREM 5.5} The canonical action of a topological group $G$
on its greatest ambit, $\Cal S(G)$, determines an embedding of
topological groups, $G\hookrightarrow\Homeo(\Cal S(G))$.
\newqed\endproclaim

The terminology we suggest using on such occasions is that the action
is {\it topologically effective.}

Presently the standard
reference to the concept of the greatest ambit is Brook's paper
\cite{Br}.

\head 6. The Universal Minimal Flow\endhead
Let $G$ be a topological group.
It is clear that a compact $G$-flow, $X$, is minimal, that is, admits
no proper compact $G$-subflows, if and only if the orbit of every
$x\in X$ is everywhere dense in $X$. Since for an arbitrary compact
$G$-flow the family of all compact non-empty 
$G$-subflows survives the 
intersections of chains, Zorn's lemma applies and yields the following result.

\proclaim{PROPOSITION 6.1} Every compact $G$-flow contains a minimal
$G$-subflow.
\newqed\endproclaim

A minimal $G$-subflow of the greatest ambit, $\Cal S(G)$, is denoted
by $\Cal M(G)$ and called the {\it universal minimal $G$-flow}.
The name is justified by the following result.

\proclaim{PROPOSITION 6.2} For a topological group $G$ every minimal
$G$-flow, $X$, is the image of the universal minimal flow, $\Cal M(G)$,
under a surjective morphism of $G$-flows.
\endproclaim

\demo{PROOF} Let $x\in X$ be any. The restriction to
(an isomorphic copy of) $\Cal M(G)$ of the unique morphism of $G$-ambits
$\Cal S(G)\to X$ sending $e$ to $x$ has a compact (therefore closed)
image in $X$, which image contains the everywhere dense
$G$-orbit of $x$ and thus coincides with $X$.
\newqed\enddemo
\smallskip

The ambiguity in choosing an $x\in X$ in the above proof
is responsible for the fact that there is, generally speaking,
no {\it canonical} morphism from $\Cal M(G)$ onto $X$, unlike it is
the case with the greatest ambit. This is also the reason why the proof
of uniqueness is slightly more involved, and we refer the reader
to e.g. \cite{A}.

\proclaim{PROPOSITION 6.3} For a topological group $G$, 
the universal minimal $G$-flow is unique up to an isomorphism.
\newqed\endproclaim

What is the `size' of the universal minimal flow for a given
topological group? In particular, when is $\Cal M(G)$ nontrivial and
when $G$ acts upon $\Cal M(G)$ effectively? By far the strongest
result in this direction to date is the following.

\proclaim{THEOREM  6.4 {\rm (Veech \cite{V1})}}
Every locally compact group $G$ acts freely on the universal
minimal flow.
\newqed\endproclaim

The Veech's proof is rather technical. While we are unable
to find a `soft' version of it, we offer a simple proof of the weaker
statement: every locally compact group acts {\it effectively} on the
universal minimal flow. This will be achieved by means of a
new description of topological groups admitting effective 
minimal flows.

\head 7. Syndetic Sets and The Size of Universal 
Minimal Flow\endhead
A subset $S$ of a topological group $G$ is called
{\it left syndetic} if for some compact $K\subseteq G$ one has
$KS=G$. If a set $S$ is left syndetic with respect to the discrete
topology on a group (that is, $FS=G$ for some
$F$ finite), then $S$ is called {\it discretely left syndetic}
or {\it big on the left}. If either $G$ is abelian or $S$ is
symmetric then `left' in both definitions can be dropped.

Here is how syndetic sets emerge in 
topological dynamics. Suppose a topological group $G$ acts on a minimal
compact $G$-flow $X$ and let $x\in V\subseteq X$, where $V$ is open.
Then the open set 
$\widetilde V=\{g\in G\colon gx\in V\}$ is big on the left
in $G$. Indeed, since the translations $hV$, $h\in G$ form 
an open cover of $X$
(assuming a $y\in X$ is not in their union, the orbit of $y$ would 
miss $V$, contrary to everywhere density in $X$), there is
a finite subcover $\{fV\colon f\in F\}$, where
$F\subseteq G$, $\vert F\vert<\infty$.
It remains to notice that $\widetilde{fV}=f\widetilde V$ for each
$f\in G$ and therefore $F\widetilde V=G$.
(To see that $\widetilde V$ is not necessarily {\it right} syndetic,
it is enough to
study the full homeomorphism group acting on the circle $\Bbb S^1$.)

\proclaim{THEOREM  7.1} Let $G$ be a topological group and $g\in G$.
The following are equivalent.
\item{\rm (i)} The motion of the universal minimal flow
$\Cal M(G)$ determined by $g$ is
nontrivial, that is, $gx\neq x$ for some $x\in\Cal M(G)$.
\item{\rm (ii)} There is a left syndetic set $S\subseteq G$ with
$g\notin\overline{SS^{-1}}$.
\item{\rm (iii)} There is a set $S\subseteq G$ big on the left
and such that 
$g\notin\overline{SS^{-1}}$.
\item{\rm (iv)} There is an open set $S\subseteq G$ big on the left and
such that
$g\notin\overline{SS^{-1}}$.
\endproclaim

\demo{PROOF} (i) $\Rightarrow$ (iv): 
Let $\Cal U_X$ denote the (unique) uniform structure
on $X=\Cal M(G)$.
For every $E\in\Cal U_X$,  
let $\tilde E=\{(g,h)\in G\times G\colon (gx,hx)\in E\}$.
The collection $\{\tilde E\colon E\in\Cal U_X\}$ forms a basis of
a uniform structure, $\Cal V$, on $G$ contained in $\Cal U_\Rsh(G)$.

Let $x\in X$ be such that $gx\neq x$. 
For some $E\in\Cal U_X$ one has $E(x)\cap E(gx)=\emptyset$.
Since $\Cal V\subseteq \Cal U_\Rsh(G)$, for
a suitable symmetric $\Cal O\in\Cal N(G)$ and an $E_1\in\Cal V$, one has 
$\Cal O_\Rsh\circ E_1\subseteq E$ and 
therefore $\Cal O E_1(x)\cap \Cal O E_1(gx)=\emptyset$.
Set $V=E_1(x)\cap g^{-1} E_1(gx)$; it is a neighbourhood of $x$ 
and $\Cal O V\cap gV\subseteq\Cal O V\cap \Cal O gV=\emptyset$.
The set $S=\widetilde V\subseteq G$ is open, big on the left, and 
$\Cal O\widetilde V=
\widetilde{\Cal OV}$, which implies that $\Cal OS\cap
gS=\emptyset$, that is, 
$g\notin\overline{SS^{-1}}$, as desired.
\smallskip

(iv) $\Rightarrow$ (iii) $\Rightarrow$ (ii): trivial.
\smallskip

(ii) $\Rightarrow$ (iv): Find a symmetric $V\in\Cal N(G)$ with 
$\overline{SS^{-1}}\cap VgV=\emptyset$ and then a $W\in\Cal N(G)$ 
such that $W^2\subseteq V$. One has $WgW\cap (WS)(WS)^{-1}=\emptyset$.
Choose a compact $K$ with $KS=G$ and then a finite $F$ with 
$FW\supseteq K$. Then $F(WS)=(FW)S\supseteq KS=G$, and the set $WS$
is open, big on the left, and $g\notin\overline{(WS)(WS)^{-1}}$.
\smallskip

(iv) $\Rightarrow$ (i):
Let $U\in\Cal N(G)$ be symmetric and such that $UgU\cap SS^{-1}=\emptyset$.
It means that $gUS\cap US=\emptyset$.
Find a right-invariant continuous pseudometric $\rho$ on $G$ with
$\rho\leq 1$ and $\Cal O_1^\rho(e)\subseteq U$.
For any $x\in G$, set $f(x)=\rho\text{-}\dist(x, G\setminus US)$.
Clearly, $0\leq f(x)\leq 1$. Whenever 
$xy^{-1}\in\Cal O_\e^\rho(e)$,
one has by the triangle inequality
$\vert f(x)-f(y)\vert\leq\rho(x,y)<\e$, that is, $f$
is in $C^\ast_\Rsh(G)$.
If $x\in S$ then $f(x)=1$, and if $x\notin US$
then $f(x)=0$.
Therefore, $f\vert_{gS}\equiv 0$ and $f\vert_{S}\equiv 1$.

Choose a finite $F\subseteq G$ with $FS=G$. 
Denote by $[\cdot]$ the closure operator in $\Cal S(G)$; since
$G=\cup_{\kappa\in F}\kappa S$ and $F$ is finite, one has 
$\Cal S(G)=\cup_{\kappa\in F}[\kappa S]$.
For each $\kappa\in F$, one clearly has 
$\overline{L_\kappa(f)}\vert_{[\kappa S]}\equiv 1$ and 
$\overline{L_\kappa(f)}\vert_{(\kappa g\kappa^{-1})[\kappa S]}\equiv 0$,
where the horizontal bar denotes the unique continuous extension of
a $\Cal U_\Rsh(G)$-u.c. bounded function to $\Cal S(G)$.

For at least one $\kappa\in F$ one has 
$[\kappa S]\cap\Cal M(G)\neq\emptyset$, where $\Cal M(G)$ is
identified with a subflow of the greatest ambit.
Choose an $x$ from the latter
intersection. The property from the previous paragraph implies that
$(\kappa g\kappa^{-1})x\neq x$ (because these two points are separated
by the function $\overline{L_\kappa(f)}$) and therefore $g(\kappa^{-1}x)\neq
\kappa^{-1}x$, where $\kappa^{-1}x\in\Cal M(G)$.
\newqed\enddemo

\head 8. Trivial Universal Minimal Flows\endhead
Let us now examine the issue of (non)trivial\-ity of the universal
minimal flow.
A topological semigroup $S$ is called {\it extremely (left) amenable}
\cite{Gr1, Mi, Gr2} if every compact $S$-flow has a
fixed point. For  topological groups this amounts to saying that $\Cal
M(G)=\{\ast\}$ is a singleton. This property is a considerably
strengthened version of amenability,
as can be seen from a well-known characterization of amenable
topological groups as those whose every continuous action upon a
convex compact set by affine transformations has a fixed point
\cite{G1}. (There are many descriptions of extremely amenable
topological semigroups, including even one in terms of geometry
of Banach spaces \cite{Gr2}.)
While examples of extremely amenable topological semigroups
abound, the problem of existence of extremely amenable topological
groups, explicitely stated by Mitchell in 1970 \cite{Mi}, happened to be
much more difficult. The first such example was published in 1974
by Herer and Christensen \cite{HR}, and
in 1996 more examples with different combinations of properties were
published
by Eli Glasner \cite{G3}, who kept his construction unpublished
for a while and observed that it was obtained independently by
Furstenberg and B. Weiss, and by the present author \cite{P}. 

Granirer proposed in 1967 \cite{Gr1} the problem of characterizing 
extremely amenable topological groups in intrinsic terms.
We suggest the following solution, which is an immediate corollary 
(and in fact is an equivalent form of) Theorem 7.1.

\proclaim{THEOREM  8.1 {\rm (Extreme Amenability Criterion)}}
A topological group $G$ is extremely amenable
if and only if whenever $S\subseteq G$ is big on the left,
$SS^{-1}$ is everywhere dense in $G$.
\newqed\endproclaim

\proclaim{COROLLARY  8.2 {\rm (Glasner \cite{G3})}}
If there exists a monothetic minimally almost periodic
topological group $G$ that is not extremely amenable,
then there exists a big subset $S\subseteq\Z$ such that
$S-S$ is not a Bohr neighbourhood of zero.
\endproclaim

\demo{PROOF} Minimal almost periodicity of a topological group
is obviously equivalent to the property:
every neighbourhood of identity in the finest precompact topology on
the underlying discrete group is everywhere dense in the original
topology on $G$.
Now it is enough to apply Theorem 8.1 to an everywhere dense
copy of $\Z$.
\newqed\enddemo
\smallskip

The question of existence of a monothetic minimally almost periodic
non-extremely amenable topological group, suggested
by Glasner, remains open, as is the 
classical question on the existence of a big set, $S$, of
integers such that $S-S$ is not a Bohr neighbourhood of zero.
(Cf. \cite{G3, V1}.)

\head 9.  Extreme Amenability of Levy Groups\endhead
The following is a slight extension of a definition from
\cite{G3}, where we drop the metrizability restriction. We say that
a topological group $G$ is a {\it Levy group} if $G$ contains an
increasing chain of compact subgroups $G_i,i\in\N$ with the
everywhere dense union and such that whenever $A_i\subseteq G_i$
have the property that $\liminf\mu_i(A_i)>0$, then for every
neighbourhood of identity, $V$, $\lim \mu_i(VA_i\cap G_i)=1$,
where $\mu_i$ denotes the normalized Haar measure on $G_i$.

An example of a Levy group (due to Glasner and, independently,
Fusternberg and Weiss)
is the group of all measurable
$\Bbb S^1$-valued complex functions on a nonatomic Lebesgue measure
space $X$ equipped with the metric
$$d(f,g)=\int \vert f(x)-g(x)\vert d\mu(x)$$
Here one can choose as
the compact subgroups $G_i$ tori of increasing dimension formed by
step functions corresponding to a sequence of refining partitions of
$X$. A subtle argument \cite{G3} shows that this group is also
monothetic.

The following was proved by Gromov and Milman \cite{GM} under the
additional condition of equicontinuity of the action, later 
removed by Glasner \cite{G3} and (independently, unpublished) 
by Furstenberg and
B. Weiss. (The assumed metrizability of the group is also
redundant.)

\proclaim{THEOREM  9.1} Every Levy group is extremely amenable.
\endproclaim

\demo{PROOF} Let $S$ be any subset big on the left, and let
$F$ be finite with $FS=G$; one can assume that $F\subseteq\cup G_i$. 
Then $\liminf\mu_i(G_i\cap S)\geq 1/n$,
where $n=\vert F\vert$. (Indeed, $\mu_i(G_i\cap S)\geq 1/n$ whenever
$F\subseteq G_i$.)  Since $G$ is a Levy group,
for any neighbourhood of identity, $V$, one has 
$\lim \mu_i(V(G_i\cap S)\cap G_i)=1$. For all $i$ large
enough, $\mu_i(((V(G_i\cap S)\cap G_i)\cdot (S^{-1}\cap G_i)))=1$, since
assuming the contrary implies that for a confinal set of $i$, a suitable
translate of $G_i\cap S$ in $G_i$, having measure $\geq 1/n$,
does not meet the set $V(G_i\cap S)\cap G_i$ of measure $>1-1/n$, 
a contradiction. 
Since $VSS^{-1}\cap G_i\supseteq
(V(G_i\cap S)\cap G_i)\cdot (S^{-1}\cap G_i))$, 
the set $VSS^{-1}\cap G_i$ has full measure and
is everywhere dense in $G_i$ starting from some $i$
and thus $VSS^{-1}$ is everywhere dense in $G$. Since $V$ is
arbitrary, $SS^{-1}$ is everywhere dense in $G$, and Theorem 8.1
finishes the proof.
\newqed\enddemo
\smallskip

\head 10. Extreme Amenability of Groups of Order Automorphisms\endhead
Let us say that a group $G$ of order automorphisms of a
linearly ordered set $X$ is $\omega${\it -transitive} if it takes
any finite subset to any other subset of the same size.

\proclaim{THEOREM 10.1 {\rm (Pestov \cite{P})}}
An $\omega$-transitive group of order automorphisms of an
infinite linearly ordered set $X$, equipped with the topology of
simple convergence on $X$, is extremely amenable.
\endproclaim

\demo{PROOF} Stabilizers, $\St_M$, 
of finite subsets $M\subseteq X$ are open
subgroups forming a neighbourhood basis at identity, and the right
factor space $G/\St_M$ can be naturally identified
with the collection $X^{(n)}$ of all ordered $n$-subsets of $X$, where 
$n=\vert M\vert$. 
Let $F\subseteq G$ be finite and such that $FS=G$. According to
Infinite Ramsey's Theorem \cite{GRS}, for some 
infinite $A\subseteq X$ and
some $f\in F$ one has $A^{(n)}\subseteq \pi_M(fS)$, where
$\pi_M\colon G\to X^{(n)}$ is the factor-map. The set
$B=f^{-1}(A)$ is also infinite and $B^{(n)}\subseteq \pi_M(S)$.
In other words, $\St_MS$ contains all transformations taking $M$ to
a subset of $B$. Consequently, $S^{-1}\St_M$ contains all
transformations such that the image of $B$ contains $M$. 
Since $B\setminus M$ is
infinite, any transformation from $G$
can be represented as a composition of
one from $S^{-1}\St_M$ and one from $\St_MS$, that is,
$\St_M SS^{-1}\St_M=G$. Since $M$ is arbitrary, it means that
$G=\overline{SS^{-1}}$, as required.
\newqed\enddemo
\smallskip

An example of such a topological group is the group
$\Aut(\Q)$ of all order automorphisms of the rational numbers, equipped
with the topology of pointwise convergence on the discrete set $\Q$. 
Since this group naturally
embeds as a
topological subgroup into the unitary group 
$\operatorname{U}(l_2(\vert\Q\vert))$
with the strong operator topology,
we conclude that the unitary group of an infinite-dimensional
Hilbert space admits no free compact flow. The same conclusion
and for the same reason holds
for the full group of permutations, $\Sym(X)$, of an infinite set $X$
equipped with the topology of simple convergence.

A further series of examples 
is obtained by observing that images of 
extremely amenable topological groups under continuous
homomorphisms onto are extremely amenable. 
This is how one deduces extreme amenability of the groups of
orientiation-preserving homeomorphisms endowed with the compact-open
topology, $\Homeo_+(I)$, $\Homeo_+\R$, and the stabilizer of
any point $s\in\Bbb S^1$ in $\Homeo_+(\Bbb S^1)$.

From the latter fact it is easy to deduce that $\Bbb S^1$ forms the
universal minimal flow for the topological
group $\Homeo_+(\Bbb S^1)$, and thus we
arrive at an example showing that a topological group can act on the
universal minimal flow effectively but not freely.

Since the free group on two generators, $F_2$, embeds into the
multiplicative group of a countable linearly ordered skewfield
\cite{N}, and any such skewfield is order isomorphic to
$\Q$, one can easily deduce that 
$F_2$ embeds into $\Aut(\Q)$ as a closed discrete topological
subgroup. It means that a closed topological subgroup
of an amenable (even extremely amenable) topological group need not
be amenable, unlike it is in the locally compact case
(cf. \cite{Gre}). 

\proclaim{CONJECTURE 10.2} Every topological group is isomorphic
to a topological subgroup of an extremely amenable topological
group.
\quest\endproclaim

All the above examples are presented in detail in \cite{P}.

\head  11. Enveloping Semigroup of a Flow and Ellis' Problem\endhead
One of the guiding problems of topological dynamics
is that of decomposing generic $G$-flows into some basic, simpler 
building blocks. 
Many important structure theorems known in topological
dynamics serve this purpose.
Here we show that a 1969 problem by Ellis on the enveloping
semigroup of the universal minimal flow is of the same kind, and
discuss its solution.

The following construction belongs to Ellis \cite{El, E3}.
Let $G$ be a topological group, and let $X$ be a compact $G$-flow.
The {\it enveloping semigroup}, $\Cal E(X)$, of $X$
is a $G$-ambit defined as follows. 
Recall that every element $g\in G$ determines a self-homeomorphism,
$\tau_g$, of $X$, called the $g${\it -motion.}
The underlying topological space
of $\Cal E(X)$ is the closure of the collection of all $g$-motions,
$g\in G$, in $X^X$ equipped with the topology of
simple convergence. Since $X^X$ is compact, so is $\Cal E(X)$.
Since every element of $\Cal E(X)$ is a (not necessarily continuous)
map $X\to X$, a natural
left action of $G$ upon $\Cal E(X)$ is determined by composing
an $f\in\Cal E(X)$ with the $g$-motion on the left:
$$\Cal E(X)\ni f\overset g\to\mapsto \tau_g\circ f\in\Cal E(X)$$
The correctness of the definition and the continuity of the 
action $G\times \Cal E(X)\to\Cal E(X)$ are checked easily.
The identity map $\operatorname{Id}_X$ has an everywhere dense orbit
in $\Cal E(X)$ (it is simply the set of all motions), and therefore
$\Cal E(X)=(\Cal E(X),\operatorname{Id}_X)$ is a $G$-ambit.
Moreover, $\Cal E(X)$ supports an obvious natural semigroup structure,
but we are not interested in it.

The following explains the significance of the
enveloping semigroup.

\proclaim{THEOREM 11.1} The enveloping semigroup $\Cal E(X)$ of
a $G$-flow $X$ is the
greatest $G$-ambit with the property that morphisms into $X$ separate
points. In other words, morphisms of $G$-flows $\Cal E(X)\to X$
separate points in $\Cal E(X)$, and whenever
$(Z,z)$ is a $G$-ambit such that
morphisms of $G$-flows $Z\to X$ separate points in $Z$, there exists
a unique morphism of $G$-ambits 
$(\Cal E(X),\operatorname{Id}_X)\to (Z,z)$.
\endproclaim

\demo{PROOF} The first statement is easy to verify: if 
$f,g\in\Cal E(X)$ and $f\neq g$, then for some $x\in X$ one must
have $f(x)\neq g(x)$ and therefore $\hat x(f)\neq \hat x(g)$,
where $\hat x\colon \Cal E(X)\to X$ is the evaluation map at $x$.
It remains to observe that $\hat x$ is continuous 
(by the very definition of the topology of simple convergence) 
and equivariant (by the definition of composition of mappings), 
that is, a morphism of $G$-flows.

Now assume $(Z,z)$ is a $G$-ambit whose points are separated by
continuous equivariant maps to $X$. Denote by $\Cal F$ the collection
of all such maps $Z\to X$. Any $f\in\Cal F$ is uniquely
determined by the image of $z$, and thus the map
$$\Cal F\ni f\mapsto f(z)\in X $$
is an injection, identifying $\Cal F$ with an everywhere dense subset of $X$.
To every $a\in Z$, associate an element of $X^\Cal F$ of the form
$f\mapsto f(a)$; the resulting mapping $Z\to X^\Cal F$ is one-to-one
and continuous and therefore (as $X$ is compact)
a homeomorphic embedding.
Finally, it is easy to verify that
the natural projection $X^X\to X^\Cal F$, dual to the inclusion
$\Cal F\hookrightarrow X$, determines a continuous equivariant map
from $\Cal E(X)\subseteq X^X$ onto $Z$. 
\newqed\enddemo
\smallskip

\definition{EXAMPLES 11.2} \ {\bf 1.}
The enveloping semigroup $\Cal E(\Cal S(G))$ of the greatest
ambit of a topological group $G$ is canonically isomorphic to the
greatest ambit $\Cal S(G)$ itself. 
\smallskip\noindent
{\bf 2.} Here is a more interesting case
where Theorem 11.1 applies. Let $G$ be a discrete group.
The {\it shift system} over $G$ is topologically
a (compact, zero-dimensional) Cantor cube $\Z_2^G$, 
upon which $G$ acts by left translations. Such $G$-flows form the 
subject of study of {\it symbolic dynamics.} 
It is shown \cite{G2, Lemma 4.1} that, rather remarkably, the
enveloping semigroup $\Cal E(\Z_2^G)$ is isomorphic to $\Cal S(G)$
and therefore the shift system in a sense carries the complete
dynamical information pertaining to a discrete group.
\tri\enddefinition

\proclaim{COROLLARY 11.3} For a topological group $G$,
the following are equivalent.
\item{\rm (i)} The canonical morphism $\Cal S(G)\to\Cal E(\Cal M(G))$
is an isomorphism of $G$-ambits.
\item{\rm (ii)} Points of $\Cal S(G)$ are separated by
continuous equivariant mappings to minimal $G$-flows.
\endproclaim

\demo{PROOF} $\Rightarrow$: follows from Theorem 11.1.
$\Leftarrow$: Let $x,y\in\Cal S(G)$ be arbitrary, and let them be
separated by a continuous equivariant 
mapping, $f$, to a minimal $G$-flow, $X$.
Fix a morphism $i\colon\Cal M(G)\to X$, which is necessarily onto,
and choose any $a\in i^{-1}(f(e))$. There exists a unique morphism of
$G$-flows $\varphi\colon \Cal S(G)\to \Cal M(G)$ sending $e$ to $a$.
Since $i(\varphi(e))=f(e)$ and therefore $i\circ\varphi=f$, one has
$\varphi(x)\neq\varphi(y)$, as required.
\newqed\enddemo

Ellis asked in 1969 \cite{E3}
whether or not the condition (i) was true for
every topological group $G$. In view of the above result, it
amounts to asking whether or not the greatest ambit of a topological group
can be reconstructed using only minimal $G$-flows, that is, is
`residually minimal.'
Now that the existence of extremely amenable groups has been
established, it is clear that in general Ellis' question 
is answered in the negative. Indeed, while the greatest ambit
$\Cal S(G)$ is always nontrivial whenever so is $G$, the universal
minimal flow of an extremely amenable topological group is a singleton
and maps to it do not separate points.

The answer to Ellis' question is positive
for every precompact group $G$, because, as we observed earlier,
in this case $\Cal M(G)$ coincides with all of $\Cal S(G)$ and is
isomorphic, as a $G$-flow, to the compact group completion of $G$ equipped
with the natural left action of $G$ by translations. Surprisingly, it
remains the only known class of topological groups answering Ellis'
question in the positive. The limitations of space prevent us
from presenting a proof of the following recent result by Glasner 
\cite{G4}, since this observation is based on deep results 
\cite{F, GW}.

\proclaim{THEOREM  11.4} For the group $\Bbb Z$ with the discrete topology,
the points of the greatest ambit $\Cal S(\Z)\cong\beta\Z$ are not
separated by continuous equivariant maps to the universal minimal flow
$\Cal M(\Z)$.
\newqed\endproclaim

In view of this result, we put forward the following.

\definition{CONJECTURE  11.5} For a topological group $G$ the following are
equivalent: (i) continuous equivariant maps to minimal $G$-flows
separate points in the greatest $G$-ambit, (ii) $G$ is precompact.
\quest\enddefinition

\head  12. Effective Minimal Flows \endhead 
The following result,
deduced at once from Theorem 7.1, is in
fact an equivalent twin of Theorem 8.1, recovered from it 
by applying the latter result to the
image of a topological group $G$ in $\Homeo(\Cal M(G))$
under the action homomorphism.

\proclaim{THEOREM  12.1}
A topological group $G$ acts effectively on its universal minimal flow
$\Cal M(G)$ if and only if for every $g\in G$, $g\neq e$, there is
a left syndetic subset $S\subseteq G$ such that
$g\notin\overline{SS^{-1}}$.
\newqed\endproclaim

As the first application of this criterion, we are able to offer a
simple proof of the following weakened form of Veech's theorem.

\proclaim{THEOREM  12.2}
A locally compact group acts effectively on its universal minimal flow.
\endproclaim

\demo{PROOF}
Let $G$ be a locally compact
group, let $g\in G$ and $g\neq e$. We will verify the condition from
Theorem 12.1. Fix a compact neighbourhood of
$g$, say $V$. Assume $V\not\ni e$. 
Zorn's Lemma implies the existence of
a maximal subset $A\subseteq G$ with $e\in A$ and 
$AA^{-1}\cap V=\emptyset$. We claim that $(V^{-1}\cup V)A=G$ and
in particular $A$ is left syndetic, which finishes the proof.

Assuming it is not so, there is some $g\in G$ with
$$g\notin (V^{-1}\cup V)A\equiv V^{-1}A\cup VA,$$
which means that $Ag^{-1}\cap V = \emptyset = gA^{-1}\cap V$.
Now we have
$$(A\cup \{g\})(A\cup \{g\})^{-1}= (AA^{-1}\cap V)\cup 
(Ag^{-1}\cap V) \cup (gA^{-1}\cap V) \cup (\{e\}\cap A) =\emptyset, $$
in contradiction with the assumed maximality of $A$. 
\newqed\enddemo
\smallskip

\definition{PROBLEM  12.3} Describe in intrinsic terms
those topological groups acting {\it freely} on their universal minimal
flows and use this description to obtain a `soft' proof of Veech's
theorem.
\quest\enddefinition

Remember that $\Homeo_+(\Bbb S^1)$ provides
an example of a topological group whose action on
the universal minimal flow is effective but not free (Section 11).

\head 13. On the Existence of Unitary Representations\endhead
Out of multitude of ways to define an amenable topological
group \cite{G1, Gre}, the one most apt for
abstract topological dynamics
is this: a topological group $G$ is amenable if and only if every
compact $G$-flow $X$ admits an invariant mean, that is,
a positive linear functional $\varphi\colon C(X)\to\C$ of norm
$1$ such that $\varphi(L_gf)=\varphi(f)$ for all $f\in C(X)$, $g\in G$.
(Cf. e.g. \cite{A}.)

 The following is a somewhat strengthened form of
Theorem 4 from \cite{HR}.

\proclaim{THEOREM 13.1} Let an amenable topological group $G$ act
effectively on the universal minimal flow $\Cal M(G)$. Then
$G$ admits a faithful strongly continuous representation in a
Hilbert space.
\endproclaim

\demo{PROOF}
Let $\varphi\colon C(\Cal M(G))\to\C$ be an invariant mean.
Set for all $f,h\in C(\Cal M(G))$: $\langle f,g\rangle=\varphi (f\overline h)$.
Clearly, $\langle,\rangle$ is an invariant
sesquilinear form on $C(\Cal M(G))$. 
Denote by $\Cal H=\widehat{C(\Cal M(G))/\Cal N_\varphi}$ the associated Hilbert
space, where $\Cal N_\varphi=\{x\in C(\Cal M(G))\colon \langle x,x\rangle=0\}$.
The group $G$ acts on $\Cal H$ by isometries. This action is strongly
continuous. Indeed, so is
the action of $G$ on $C(\Cal M(G))$ and since
the pre-Hilbert topology on $C(\Cal M(G))$ is coarser than
the uniform topology, each orbit map 
$G\ni g\mapsto g\cdot f\in\Cal H$ is continuous
whenever $f\in C(\Cal M(G))$. It means that the representation of $G$ in
$\Cal H$ is continuous with respect to the topology of simple convergence
on an everywhere dense subset of $\Cal H$. But on the unitary
group $U(\Cal H)$ any such
topology coincides with the strong operator topology.

It remains to show that the representation is faithful. Let $g\in G$ be
any, and let $x\in\Cal M(G)$ be such that $gx\neq x$. 
Fix a neighbourhood $V\ni x$ with $gV\cap V=\emptyset$. There exists
an $f\in C(\Cal M(G))$ with $0\leq f\leq 1$, $\supp f\subset V$, and
$\varphi(f)>0$.
Indeed, otherwise the support of the invariant mean $\varphi$ would not
meet any translate $gV$, and since finitely many of those cover
$\Cal M(G)$ by force of its minimality, the support of $\varphi$ would
be empty, in contradiction to $\norm\varphi=1$.
Now it is easy to see that
$\norm{f-g\cdot f}=2\norm{f}>0$, which means that
$g$ as an operator in $\Cal H$ is different from identity.
\newqed\enddemo

Inverting the above result, one concludes that an amenable topological
group admitting no nontrivial strongly continuous unitary representations
is extremely amenable. This observation (made in the particular case of
abelian topological groups, which are all known to be amenable)
was employed in \cite{HR} by
Herer and Christensen, who constructed an abelian topological group
without nontrivial unitary representations and deduced its extreme
amenability.

Combining Theorems 13.1 and 7.1 leads to the following new
existence result for unitary representations.

\proclaim{COROLLARY 13.2} Let $G$ be an amenable topological group
and let $g\in G$ be such that for some left syndetic subset
$S\subseteq G$ one has $g\notin\overline{SS^{-1}}$. Then $g$ is
separated from identity by a strongly continuous unitary representation
of $G$.
\endproclaim

\demo{PROOF} Denote the action of $G$ on $\Cal M(G)$ by $\tau$.
According to Theorem 7.1, $\tau(g)\neq e$ in the 
topological group $\tau(G)<\Homeo(\Cal M(G))$. Since the latter group
acts effectively and minimally on $\Cal M(G)$, it admits a faithful
unitary representation $\pi$ by 13.1, and the unitary representation
$\pi\circ\tau$ of $G$ is strongly continuous and separates $g$ from
identity.
\newqed\enddemo

Unfortunately, the above result falls short of a criterion, since
even for amenable topological groups the existence of
such a left syndetic set $S$ is not
necessary for the existence of a nontrivial unitary representation.
Indeed, as observed in \cite{G4},
the monothetic extremely amenable group from Glasner--Furstenberg--Weiss
example (Section 10) embeds as a
topological subgroup into the unitary group of $l_2$ with the strong
operator topology. Thus, the old problem of describing in
intrinsic terms those topological groups admitting a complete system
of strongly continuous unitary representations remains open.

In view of Herer--Christensen's example \cite{HC},
the following curious question seems to be unanswered.

\definition{QUESTION 13.3} 
Does there exist a {\it monothetic} topological group admitting no
nontrivial strongly continuous unitary representations in a
Hilbert space?
\quest\enddefinition

The importance of big sets in abstract harmonic analysis and
representation theory is underlined by the following classical result.

\proclaim{THEOREM  13.4 {\rm (Cotlar and Ricabarra \cite{CR})}}
Let $G$ be an abelian topological group. An element $g\in G$ is 
separated from identity by a continuous character
if and only if there exists a big symmetric
open set $S\subseteq G$ with $g\notin S^6$.
\newqed\endproclaim

Later Ellis and Keynes \cite{EK} reduced the number $6$ to $4$.

In this connection, the following might be interesting.

\proclaim{CONJECTURE 13.5} 
Let $G$ be a topological group. An element $g\in G$ is 
separated from identity by a continuous finite-dimensional
unitary representation if and only if there exists a big symmetric 
open set $S\subseteq G$ with $g\notin S^4$ (or $S^6$, etc.).
\quest\endproclaim

Of course, the necessity ($\Rightarrow$) is always
valid. Possibly, one extra condition to be imposed on $S$
is that of being {\it invariant} under inner automorphisms of $G$.

\enddocument
\bye